\documentclass[12pt]{article}
\def\C{\mathbf{C}}
\def\const{\mathrm{const}}
\def\grad{\mathrm{grad}\ }
\def\z{\mathbf{z}}
\def\w{\mathbf{w}}
\begin{document}
\title{Tangencies between holomorphic maps and holomorphic
laminations}
\author{A. Eremenko\footnote{Supported by NSF grant
DMS-0555279.}$\;$
and A. Gabrielov\footnote{Supported by the NSF grant
DMS-0801050.}}
\date{\today}
\maketitle
\begin{abstract}
We prove that the set of leaves of a holomorphic
lamination of codimension one that are tangent to a germ
of a holomorphic map is discrete.
\end{abstract}

Let $F$ be a holomorphic lamination of codimension
one in an open set $V$
in a complex Banach space $B$. In this paper, this means
that $V=W\times\C$,
where $W$ is a neighborhood of the
origin in some Banach space,
and the leaves $L_\lambda$ of the
lamination are disjoint graphs of holomorphic
functions $\w\mapsto f(\lambda,\w), \; W\to\C$.
For holomorphic functions in a Banach space
we refer to \cite{H}.
Here $\lambda$ is a parameter and we assume
that the dependence of $f$ on $\lambda$ is continuous.
A natural choice of this parameter is such that
$\lambda=f(\lambda,0)$, in which case the continuity 
with respect to $\lambda$ follows from the so-called
$\lambda$-lemma of Mane-Sullivan-Sad and Lyubich,
see, for example \cite{H}. With this choice of the parameter,
our definition of a lamination coincides with
that of a holomorphic motion of $\C$ parametrized by $W$.

Let $\gamma:U\to V$ be a holomorphic map, $U\subset\C^n$. 
We say that $\gamma$ is tangent
to the lamination at a point $\z_0\in U$ if
the image of the derivative $\gamma'(\z_0)$
is contained in the tangent
space $T_L(\gamma(\z_0)),$ where $L$ is the leaf passing
through $\gamma(\z_0)$. A leaf for which this holds is
called a tangent leaf to $\gamma$. 
\vspace{.1in}

\noindent
{\bf Theorem.} {\em Let $K$ be a compact
subset of $U$. Then the set of leaves
tangent to $\gamma$ at the points of $K$
is finite.}
\vspace{.1in}

For the case of holomorphic curves ($n=1$) this result is
contained in \cite[Lemma 9.1]{0} where it is credited to Douady.
Artur Avila, in a conversation with the authors, proposed
to extend this result to
arbitrary holomorphic maps. According
to Avila, this generalization has several 
applications to holomorphic dynamics.
\vspace{.1in}

{\em Proof}. We assume without loss of generality
that $$f(0,\w)\equiv 0,$$ and
that $L_0$ is tangent to $\gamma$ at $\z_0=0$. 

We have to show that
other tangent leaves cannot accumulate to $L_0$.
Suppose the contrary, that is suppose
that there is a sequence $\lambda_k\to 0$ such that
$L_{\lambda_k}$ are tangent to $\gamma$, and
let $L_{\lambda_k}$ be the graphs of the functions 
$f_k(\w)=f(\lambda_k,\w).$ We may assume that tangency
points $\z_k\to 0$.

We make several preliminary reductions.
\vspace{.1in}

1. Let $$\gamma(\z)=
(\phi(\z),\psi(\z))\in W\times\C.$$
Consider the new lamination in $U\times\C$ whose leaves
are the graphs of $f^*(\lambda,\z)=f(\lambda,\phi(\z))$
and the new map $\gamma^*(\z)=(\z,\psi(\z))$.
Then $\gamma^*$ is tangent to a leaf $L^*$ if and only if
$\gamma$ is tangent to $L$. This reduces our problem to the
case that $W$ is an open set in $\C^n$ and the map $\gamma$ is a graph of
a function $\psi$ of the same variable as the functions 
$f_k$. From now on we assume that $U=W$ and
$\gamma(\w)=(\w,\psi(\w)).$
\vspace{.1in}

2. Now we reduce the problem to the case that
$\psi$ is a monomial. For this we use the
desingularization theorem of Hironaka \cite{1,2,3,4}.

Let $X$ be a complex analytic manifold, and $\psi$ 
an analytic function on $X$. Then there exists a
complex analytic manifold $M$ and a proper surjective map 
$\pi:M\to X$ such that the restriction of $\pi$ onto
the complement of the $\pi$-preimage of the set $\{\psi=0,\ \psi'=0\}$
is injective and for each point
$\z_0\in\pi^{-1}(\{\psi=0\})$
there is a local coordinate system with the origin
at $\z_0$ such that
$\psi\circ\pi$ is a monomial $z_1^{m_1}\ldots z_n^{m_n}$.

Let $Y=W\times\C$, and let $S\subset Y$
be the set of points $(\w,t)$ with $t\neq 0$
such that the graph of $t=\psi(\w)$ is tangent to the
lamination. In our proof by contradiction,
we assume that the origin belongs to the
closure of $S$. Let $N=M\times\C$, and let $\rho:N\to Y$
be the map defined by $\rho(\z,t)=(\pi(\z),t)$.
Then $\rho^{-1}(F)$ is the lamination whose leaves are
the components of the $\rho$-preimages of the leaves of
$F$, and the set $T=\rho^{-1}(S)$ has a limit point
$(\z_0,0)$ with $\pi(\z_0)=0$ since $\pi$ is proper.
Also the set $T$ is exactly the set of those points
in $N$ where the graph $t=\psi\circ\pi(\z)$ is tangent 
to the lamination $\rho^{-1}(F)$ since $\rho$ is injective
in a neighborhood of each point of $T$.
(Any point $(\z,t)$ where $\rho$ is not injective
satisfies $\psi(\pi(\z))=0$ while at every point of $T$
we have $\psi(\pi(\z))\neq 0$.)
This reduces our problem to the case
that $\psi$ is a monomial.
\vspace{.1in}

3. We may assume now that
$W=\{ \z:|\z|<2\}$. So we are in the following situation.
$$\psi(z)=z_1^{m_1}\ldots z_n^{m_n},$$
and $\{ f_k\}$ is a family of holomorphic functions on $W$
with disjoint graphs, $f_k(\z)\neq 0$ for $\z\in W$,
and $f_k\to 0$ as $k\to \infty$ uniformly on $W$. 
Moreover, for some sequence $\z_k\to 0$
we have 
\begin{equation}
\label{1}
f_k(\z_k)=\prod_{j=1}^n z_{j,k}^{m_j},
\end{equation}
\begin{equation}
\label{2}
\grad f_k(\z_k)=(m_1z_{1,k}^{m_1-1}z_{2,k}^{m_2}\ldots ,\;
m_2z_1^{m_1}z_{2,k}^{m_2-1}\ldots,\;\ldots)
\end{equation}
assuming zero values for the components with $m_j=0$.
We may assume that $| f_k|\leq 1$ in $W$.
Setting $f_k=\exp g_k$ we obtain that $\Re g_k\leq 0$.
Now we put
$$h_k(\z)=g_k(\z+\z_k)-g_k(\z_k).$$
Then the $h_k$ are defined in the unit ball and
satisfy
$$\Re h_k(\z)\leq\sum_{j:m_j>0} m_j\log|z_{j,k}|^{-1}.$$
From this we conclude that
\begin{equation}
\label{3}
\left|\frac{\partial h_k}{\partial z_j}(0)\right|\leq
2\,\sum_{j: m_j>0} m_j\log|z_{j,k}|^{-1}.
\end{equation}

This follows from the 
\vspace{.1in}

\noindent
{\bf Lemma}. {\em Let $h$ be an analytic function
in the unit disc, $h(0)=0$, and $\Re h\leq A$, where $A>0$.
Then $|h'(0)|\leq 2A.$}
\vspace{.1in}

This is an immediate consequence of the Schwarz Lemma.
\vspace{.1in}

On the other hand, (\ref{1}) and (\ref{2}) imply that, for $m_j>0$,
\begin{equation}
\label{4}
\left|\frac{\partial h_k}{\partial z_j}(0)\right|=
\frac{m_j}{|z_{j,k}|}.
\end{equation}
Assume without loss of generality that $m_1>0$ and
$$|z_{1,k}|=\min_{j: m_j>0}|z_{j,k}|.$$
Then the RHS of (\ref{3}) is at most
$$\const\log|z_{1,k}|^{-1},$$
while the RHS of (\ref{4}) is 
$$\frac{m_1}{|z_{1,k}|}.$$
As $|z_{1,k}|\to 0$, we obtain a contradiction
which proves our theorem.

\vspace{.2in}

{\em Purdue University

West Lafayette, IN 47907 USA

eremenko@math.purdue.edu

agabriel@math.purdue.edu}

\end{document}